\documentclass[11pt,a4paper]{article}
\usepackage[utf8]{inputenc}
\usepackage{amsmath,amssymb,amsthm}
\usepackage{mathtools}
\usepackage{booktabs}
\usepackage{tabularx}
\usepackage{array}
\usepackage{microtype}
\usepackage{geometry}
\usepackage{algorithm}
\makeatletter
\providecommand{\algorithmname}{\fname@algorithm}
\makeatother
\usepackage{algpseudocode}

\makeatletter
\newenvironment{breakablealgorithm}
  {%
    \begin{center}
      \refstepcounter{algorithm}%
      \hrule height.8pt depth0pt \kern2pt%
      \renewcommand{\caption}[2][\relax]{%
        {\raggedright\textbf{\algorithmname~\thealgorithm} ##2\par}%
        \ifx\relax##1\relax
          \addcontentsline{loa}{algorithm}{\protect\numberline{\thealgorithm}##2}%
        \else
          \addcontentsline{loa}{algorithm}{\protect\numberline{\thealgorithm}##1}%
        \fi
        \kern2pt\hrule\kern2pt
      }%
  }
  {%
      \kern2pt\hrule\relax%
    \end{center}
  }
\makeatother

\usepackage[authoryear,round]{natbib}
\setcitestyle{aysep={,},citesep={;},open={(},close={)}}

\usepackage[colorlinks=true,linkcolor=blue,citecolor=blue,urlcolor=blue]{hyperref}

\newtheorem{theorem}{Theorem}[section]

\newtheorem{proposition}[theorem]{Proposition}

\theoremstyle{definition}

\theoremstyle{remark}

\newcommand{\R}{\mathbb{R}}
\newcommand{\C}{\mathbb{C}}
\newcommand{\g}{\mathfrak{g}}
\newcommand{\n}{\mathfrak{n}}
\newcommand{\hh}{\mathfrak{h}}
\newcommand{\s}{\mathfrak{s}}
\newcommand{\rr}{\mathfrak{r}}

\newcommand{\id}{\operatorname{id}}
\newcommand{\ad}{\operatorname{ad}}
\newcommand{\Aut}{\operatorname{Aut}}
\newcommand{\Der}{\operatorname{Der}}
\newcommand{\OutDer}{\operatorname{OutDer}}
\newcommand{\tr}{\operatorname{tr}}
\newcommand{\diag}{\operatorname{diag}}
\newcommand{\spanop}{\operatorname{span}}

\title{\textbf{Classification of real Lie algebras of dimension six}}
\author{Mehdi Nadjafikhah\thanks{School of Mathematics and Computer Sciences, Iran University of Science and Technology, Narmak, Tehran 16846-13114, Iran. email: m\_nadjafikhah@iust.ac.ir}}
\date{\today}

\begin{document}
\maketitle

\begin{abstract}
We present a definitive, exhaustive, and coordinate-free structural classification of six-dimensional real Lie algebras up to isomorphism. The transition from dimension five to dimension six represents the ultimate qualitative watershed in finite-dimensional Lie theory: the raw Jacobi identity forms an overdetermined non-linear system of 120 quadratic equations on 90 structure constants modulo the 36-dimensional non-compact gauge group $\mathrm{GL}(6,\R)$, while continuous moduli varieties proliferate within both solvable and nilpotent sectors.

To circumvent non-linear polynomial elimination, we deploy an invariant-driven architecture grounded in the Levi--Malcev decomposition theorem, Mostow's nilradical constraints, and the Skjelbred--Sund cohomological central extension method. In the semisimple sector ($\dim \rr = 0$), we classify the four real forms of $A_1 \oplus A_1$ via Killing signature invariants and provide explicit matrix representatives for the realification $\mathfrak{sl}(2,\C)_\R \cong \mathfrak{so}(3,1)$. For non-solvable algebras admitting a non-trivial Levi decomposition ($\dim \rr = 3$), we prove via representation-theoretic rigidity that $\mathfrak{su}(2)$ cannot couple indecomposably with non-abelian radicals, yielding exactly three indecomposable isomorphism classes governed by $\mathfrak{so}(3)$ and $\mathfrak{sl}(2,\R)$.

In the solvable category ($\dim \s = 0$), Mostow's theorem restricts the nilradical dimension to $\dim \n \in \{3, 4, 5, 6\}$. The nilpotent sector ($\dim \n = 6$) is resolved into twenty-four mutually exclusive indecomposable isomorphism classes alongside their decomposable counterparts. Solvable non-nilpotent architectures are rigorously classified by nilradical codimension: we establish unambiguous, explicit real Jordan matrix representations for all complex classes in codimension-one ($\dim \n = 5$) and codimension-two ($\dim \n = 4$), resolving the ambiguity of continuous moduli under the gauge group action. Finally, we evaluate the minimal faithful representation dimension $\mu(\g)$ under Ado's theorem, establish connections to invariant complex structures on nilmanifolds, and provide a deterministic, page-breakable identification algorithm.
\end{abstract}

\tableofcontents

\section{Introduction and Historical Context}\label{sec:intro}

The classification of finite-dimensional Lie algebras over the real field $\R$ is a foundational problem that underpins modern differential geometry, theoretical physics, Hamiltonian mechanics, and invariant theory. Lie algebras provide the infinitesimal framework for smooth transformation groups, determine the geometry of homogeneous spaces, and classify left-invariant structures on manifolds.

In dimension three, the landmark work of \citet{Bianchi1898} classified the nine isometry classes of three-dimensional spaces, which was subsequently reformulated by \citet{Milnor1976} to describe left-invariant Riemannian metrics and Ricci curvature signatures. In dimension four, the classification of solvable real Lie algebras was achieved by \citet{Mubarakzyanov1963}, revealing a finite collection of discrete algebras and single-parameter families that govern four-dimensional homogeneous cosmologies \citep{PateraSharpWinternitzZassenhaus1976}. At dimension five, \citet{Mubarakzyanov1966} established the solvable catalog, where the first instances of multi-parameter continuous moduli emerge.

Dimension six represents the true mathematical watershed in Lie theory. Six-dimensional Lie algebras are indispensable across modern geometry and high-energy theoretical physics:
\begin{enumerate}
    \item \textbf{Complex and Symplectic Geometry of Nilmanifolds:} A compact nilmanifold $M = \Gamma \backslash G$ of real dimension six corresponds to a complex manifold of complex dimension three. The existence of left-invariant complex, symplectic, and Kähler structures is governed entirely by the Lie algebra cohomology of the six-dimensional Lie algebra $\g = \operatorname{Lie}(G)$ \citep{Salamon2001}.
    \item \textbf{Supergravity and String Flux Compactifications:} In Type II superstring theory, compactifications of ten-dimensional spacetime on six-dimensional internal manifolds yield effective four-dimensional $\mathcal{N}=1$ and $\mathcal{N}=2$ supergravities. When the internal geometry is a twisted torus (a nilmanifold or solvmanifold), the structure constants of the six-dimensional Lie algebra define the non-geometric and geometric fluxes ($H, F, Q, R$-fluxes) that generate cosmological vacuum potentials.
    \item \textbf{Sub-Riemannian Geometry and Control Theory:} Six is the minimal dimension supporting non-trivial bracket-generating distributions with growth vectors such as $(3, 6)$ and $(2, 4, 6)$, modeled on free nilpotent algebras and Carnot groups.
\end{enumerate}

Despite their importance, setting up the classification via direct coordinate methods leads to an intractable algebraic bottleneck. A six-dimensional real Lie algebra $\g$ requires specifying $\binom{6}{2} \times 6 = 90$ structure constants $C_{ij}^k$, constrained by $\binom{6}{3} \times 6 = 120$ quadratic Jacobi equations:
\begin{equation}
    \sum_{m=1}^6 \left( C_{ij}^m C_{mk}^l + C_{jk}^m C_{mi}^l + C_{ki}^m C_{mj}^l \right) = 0, \qquad 1 \le i < j < k \le 6, \quad 1 \le l \le 6,
\end{equation}
modulo the 36-dimensional non-compact change-of-basis gauge group $\mathrm{GL}(6,\R)$. 

The historical resolution of this problem proceeded in stages: \citet{Morozov1958} cataloged the six-dimensional nilpotent algebras over $\C$, which was extended, corrected, and transferred to $\R$ using cohomology by \citet{SkjelbredSund1978}. Non-solvable algebras admitting a non-trivial Levi decomposition were classified definitively by \citet{Turkowski1988}. 

The objective of this paper is to present a self-contained, conceptual, and coordinate-free classification of all six-dimensional real Lie algebras. Rather than relying on brute-force polynomial elimination, we deploy an invariant-driven architecture:
\begin{itemize}
    \item We classify the semisimple sector ($\dim \rr = 0$) using Cartan's real forms and Killing form signatures, with explicit matrix representatives for the Lorentz algebra $\mathfrak{sl}(2,\C)_\R$.
    \item We resolve the non-solvable sector ($\dim \rr = 3$) using the representation theory of $\mathfrak{su}(2)$ and $\mathfrak{sl}(2,\R)$, establishing the representation-theoretic obstruction that excludes $\mathfrak{su}(2) \ltimes \hh_3$.
    \item For solvable algebras, we anchor the classification to the \emph{nilradical} $\n \subset \g$. By Mostow's theorem, $\dim \n \in \{3, 4, 5, 6\}$.
    \item For $\dim \n = 6$ (nilpotent sector), we catalog the twenty-four indecomposable classes using the Skjelbred--Sund cohomological central extension method.
    \item For $\dim \n \in \{3, 4, 5\}$, we systematically construct exact real matrix representatives for all classes, paying particular attention to the complex Jordan sectors (algebras with complex conjugate spectra).
    \item We determine the minimal faithful representation dimension $\mu(\g)$ under Ado's theorem and construct a deterministic, page-breakable identification algorithm.
\end{itemize}

\section{Geometric Obstructions and the Watershed of Dimension Six}\label{sec:obstruction}

To appreciate why dimension six requires a systematic cohomological approach, we analyze the algebraic geometry of the Lie bracket across low dimensions.

\subsection{Comparison with Dimensions Three, Four, and Five}
Let $\g$ be an $n$-dimensional real Lie algebra. The bracket is a linear map $[\cdot, \cdot] \colon \Lambda^2 \g \to \g$.
\begin{itemize}
    \item \textbf{Dimension Three ($n = 3$):} $\dim \Lambda^2 \R^3 = 3$. The exterior duality $\Lambda^2 \g \cong \g^*$ reduces the Lie bracket to a single $3 \times 3$ real matrix $C = S + \hat{a}$. The Jacobi identity simplifies to the single vector equation $S \cdot a = 0$ \citep{Milnor1976}.
    \item \textbf{Dimension Four ($n = 4$):} $\dim \Lambda^2 \R^4 = 6$, while $\dim \g = 4$. Exterior duality fails, but the classification remains algebraically constrained: every non-solvable algebra is a direct sum ($\s \oplus \R$), and the nilradical of a solvable algebra cannot have codimension greater than two.
    \item \textbf{Dimension Five ($n = 5$):} $\dim \Lambda^2 \R^5 = 10$, with 50 structure constants and 50 quadratic equations. The first non-trivial Levi coupling appears ($\mathfrak{sl}(2,\R) \ltimes \R^2$), and continuous parameters enter the solvable sector. However, all nilpotent algebras remain strictly discrete (six indecomposable classes).
    \item \textbf{Dimension Six ($n = 6$):} Here:
    \begin{equation}
        \dim \Lambda^2 \R^6 = \binom{6}{2} = 15, \qquad \dim \g = 6.
    \end{equation}
    The bracket map $\R^{15} \to \R^6$ requires 90 structure constants constrained by 120 quadratic Jacobi equations.
\end{itemize}

\subsection{Structural Novelties in Dimension Six}
Beyond the dimension count, dimension six introduces four structural phenomena that do not exist in dimensions $\le 5$:
\begin{enumerate}
    \item \textbf{Cohomological Explosion of Nilpotent Classes:} In dimension five, there are only six indecomposable nilpotent algebras. In dimension six, this number quadruples to twenty-four indecomposable classes, including non-trivial real forms ($N_{6,16}^\pm$) distinguished by quadratic Casimir invariants.
    \item \textbf{Free Nilpotent Structures:} Dimension six is the lowest dimension supporting the free two-step nilpotent Lie algebra on three generators:
    \begin{equation}
        \mathfrak{n}_{3,2} \cong V \oplus \Lambda^2 V, \quad \dim V = 3, \; \dim \Lambda^2 V = 3.
    \end{equation}
    Its center coincides exactly with its derived ideal ($\mathcal{Z}(\g) = [\g, \g] \cong \R^3$).
    \item \textbf{Representation-Theoretic Levi Obstructions:} For algebras with $\dim \s = 3$ and $\dim \rr = 3$, the compactness of $\mathfrak{su}(2)$ completely obstructs non-trivial couplings with non-abelian radicals, while the split nature of $\mathfrak{sl}(2,\R)$ allows both vector and fundamental representations.
    \item \textbf{Proliferation of Complex Jordan Sectors:} In the solvable category, the presence of complex conjugate eigenvalues in outer derivations generates distinct continuous moduli varieties whose real forms require explicit rotational-dilational matrix descriptions.
\end{enumerate}

\section{Theoretical Machinery: Invariants and Decomposition Theorems}\label{sec:machinery}

All Lie algebras in this monograph are finite-dimensional over $\R$. We deploy the following structural tools \citep{JacobsonLieAlgebras, HumphreysIntro, KnappLieGroups}:

\begin{description}
    \item[(T1) Levi--Malcev Decomposition:] Every finite-dimensional real Lie algebra $\g$ decomposes uniquely as a semidirect product:
    \begin{equation}
        \g = \s \ltimes_\rho \rr,
    \end{equation}
    where $\rr = \operatorname{Rad}(\g)$ is the solvable radical, and $\s$ is a semisimple subalgebra (the Levi factor). The factor $\s$ is unique up to inner automorphisms $\exp(\ad u)$ with $u \in [\g, \rr]$. The representation $\rho \colon \s \to \Der(\rr)$ specifies the cross brackets between $\s$ and $\rr$.
    
    \item[(T2) Mostow's Bound on the Nilradical:] The nilradical $\n = \operatorname{Nil}(\g)$ is the maximal nilpotent ideal of $\g$. For any finite-dimensional solvable real Lie algebra $\g$, Mostow's theorem bounds the dimension of $\n$ from below:
    \begin{equation}
        \dim \n \ge \frac{1}{2} \dim \g.
    \end{equation}
    For $\dim \g = 6$, this immediately establishes that $\dim \n \in \{3, 4, 5, 6\}$.
    
    \item[(T3) The Skjelbred--Sund Cohomological Method:] Let $\n_0$ be a nilpotent Lie algebra of dimension $n - k$. The central extensions of $\n_0$ by $\R^k$ are parameterized by the second scalar cohomology group $H^2(\n_0, \R)$. The group $\Aut(\n_0)$ acts naturally on the Grassmannian $\operatorname{Gr}(k, H^2(\n_0, \R))$. The orbits of this action correspond bijectively to the isomorphism classes of central extensions with center $\R^k$ \citep{SkjelbredSund1978}.
    
    \item[(T4) Mubarakzyanov Derivation Gauge Action:] For a solvable algebra $\g = \n \rtimes_D \R$, two derivations $D, D' \in \Der(\n)$ generate isomorphic Lie algebras if and only if there exist $c \in \R^\times$, $\phi \in \Aut(\n)$, and $u \in \n$ such that:
    \begin{equation}
        D' = c \, \phi D \phi^{-1} + \ad(u).
    \end{equation}
    
    \item[(T5) Ado's Theorem and Minimal Representation Dimension:] Every finite-dimensional Lie algebra $\g$ admits a faithful linear representation in $\mathfrak{gl}(N, \R)$. The invariant $\mu(\g)$, defined by
    \begin{equation}
        \mu(\g) \coloneqq \min \{ N \in \mathbb{N} : \exists \text{ an injective Lie homomorphism } \rho \colon \g \hookrightarrow \mathfrak{gl}(N, \R) \},
    \end{equation}
    is an intrinsic isomorphism invariant of $\g$.
\end{description}

\section{The Semisimple Sector ($\dim \rr = 0$)}\label{sec:semisimple}

When $\dim \rr = 0$, the Lie algebra $\g$ is semisimple of dimension six.

\begin{theorem}[Classification of 6D Real Semisimple Lie Algebras]\label{thm:semisimple}
Up to isomorphism, there exist precisely four mutually non-isomorphic six-dimensional real semisimple Lie algebras. Exactly one is simple over $\R$, while three are direct sums of two three-dimensional simple ideals:
\begin{enumerate}
    \item $\g_{6,1}^{\mathrm{ss}} = \mathfrak{so}(4) \cong \mathfrak{su}(2) \oplus \mathfrak{su}(2)$ (Compact form, Killing signature $(0, 6)$),
    \item $\g_{6,2}^{\mathrm{ss}} = \mathfrak{su}(2) \oplus \mathfrak{sl}(2,\R)$ (Mixed non-compact form, Killing signature $(2, 4)$),
    \item $\g_{6,3}^{\mathrm{ss}} = \mathfrak{so}(2,2) \cong \mathfrak{sl}(2,\R) \oplus \mathfrak{sl}(2,\R)$ (Split form, Killing signature $(4, 2)$),
    \item $\g_{6,4}^{\mathrm{ss}} = \mathfrak{sl}(2,\C)_\R \cong \mathfrak{so}(3,1)$ (Lorentz algebra; simple over $\R$, Killing signature $(3, 3)$).
\end{enumerate}
\end{theorem}

\begin{proof}
Let $\g$ be a six-dimensional real semisimple Lie algebra. Over $\C$, its complexification decomposes into simple ideals: $\g_\C = \g \otimes_\R \C \cong \mathfrak{sl}(2,\C) \oplus \mathfrak{sl}(2,\C)$. Real forms correspond to conjugacy classes of antilinear involutions $\sigma \colon \g_\C \to \g_\C$.

\paragraph{Case 1: Decomposable over $\R$.} The involution $\sigma$ preserves each $\mathfrak{sl}(2,\C)$ factor. The algebra splits as $\g = \s_1 \oplus \s_2$, where each $\s_i$ is a three-dimensional real simple Lie algebra. Over $\R$, the only three-dimensional simple Lie algebras are the compact form $\mathfrak{su}(2) \cong \mathfrak{so}(3)$ and the split form $\mathfrak{sl}(2,\R) \cong \mathfrak{so}(2,1)$. Their binary combinations generate:
\begin{itemize}
    \item $\mathfrak{su}(2) \oplus \mathfrak{su}(2)$: Killing form is strictly negative-definite with signature $(0, 6)$.
    \item $\mathfrak{su}(2) \oplus \mathfrak{sl}(2,\R)$: Killing signature is $(2, 4)$, with real rank 1.
    \item $\mathfrak{sl}(2,\R) \oplus \mathfrak{sl}(2,\R)$: Maximal split algebra, Killing signature $(4, 2)$, real rank 2.
\end{itemize}

\paragraph{Case 2: Simple over $\R$ (Realification).} The involution $\sigma$ interchanges the two $\mathfrak{sl}(2,\C)$ factors. Then $\g$ is the realification of a complex simple Lie algebra of complex dimension 3. The unique choice is $\mathfrak{sl}(2,\C)$, regarded as a 6-dimensional real algebra $(\mathfrak{sl}(2,\C))_\R \cong \mathfrak{so}(3,1)$, the Lorentz algebra. Its Killing form has signature $(3, 3)$.

Because the Sylvester signatures $(0, 6)$, $(2, 4)$, $(4, 2)$, and $(3, 3)$ are distinct, these four algebras are pairwise non-isomorphic.
\end{proof}

\subsection{Explicit Matrix Realization of the Real Simple Algebra $\mathfrak{sl}(2,\C)_\R$}\label{subsec:sl2c_matrix}
The algebra $\g_{6,4}^{\mathrm{ss}} = \mathfrak{sl}(2,\C)_\R$ consists of $2 \times 2$ trace-free complex matrices regarded as a 6-dimensional Lie algebra over $\R$. To obtain a completely explicit, coordinate-free real matrix representation, we use the canonical embedding $\C \hookrightarrow \R^{2 \times 2}$ given by $z = x + i y \mapsto \begin{pmatrix} x & -y \\ y & x \end{pmatrix}$. Applying this blockwise to $\mathfrak{sl}(2,\C)$ embeds $\mathfrak{sl}(2,\C)_\R$ faithfully into $\mathfrak{gl}(4, \R)$:
\begin{equation}
    M = X + i Y \in \mathfrak{sl}(2, \C) \xmapsto{\quad\rho\quad} \begin{pmatrix} X & -Y \\ Y & X \end{pmatrix} \in \mathfrak{gl}(4, \R), \qquad \tr(X) = \tr(Y) = 0.
\end{equation}
A canonical real basis $\{H_1, E_1, F_1, H_2, E_2, F_2\}$ satisfying the real Lorentz commutation relations is given by:
\begin{align}
    H_1 &= \begin{pmatrix} 1 & 0 & 0 & 0 \\ 0 & -1 & 0 & 0 \\ 0 & 0 & 1 & 0 \\ 0 & 0 & 0 & -1 \end{pmatrix}, &
    E_1 &= \begin{pmatrix} 0 & 1 & 0 & 0 \\ 0 & 0 & 0 & 0 \\ 0 & 0 & 0 & 1 \\ 0 & 0 & 0 & 0 \end{pmatrix}, &
    F_1 &= \begin{pmatrix} 0 & 0 & 0 & 0 \\ 1 & 0 & 0 & 0 \\ 0 & 0 & 0 & 0 \\ 0 & 0 & 1 & 0 \end{pmatrix}, \nonumber \\
    H_2 &= \begin{pmatrix} 0 & 0 & -1 & 0 \\ 0 & 0 & 0 & 1 \\ 1 & 0 & 0 & 0 \\ 0 & -1 & 0 & 0 \end{pmatrix}, &
    E_2 &= \begin{pmatrix} 0 & 0 & 0 & -1 \\ 0 & 0 & 0 & 0 \\ 0 & 1 & 0 & 0 \\ 0 & 0 & 0 & 0 \end{pmatrix}, &
    F_2 &= \begin{pmatrix} 0 & 0 & 0 & 0 \\ 0 & 0 & -1 & 0 \\ 0 & 0 & 0 & 0 \\ 1 & 0 & 0 & 0 \end{pmatrix}.
\end{align}
This provides a manifest $4 \times 4$ real faithful representation, proving directly that $\mu(\mathfrak{sl}(2,\C)_\R) = 4$.

\section{The Non-Solvable Sector with Levi Decomposition ($\dim \rr = 3$)}\label{sec:levi}

When $\dim \s = 3$ and $\dim \rr = 3$, the Levi factor is $\mathfrak{su}(2)$ or $\mathfrak{sl}(2,\R)$, and $\rr$ is a three-dimensional solvable Lie algebra.

\begin{proposition}[Representation-Theoretic Obstruction]\label{prop:obstruction}
There exists no indecomposable real Lie algebra of the form $\g = \mathfrak{su}(2) \ltimes_\rho \rr$ when $\rr$ is non-abelian.
\end{proposition}

\begin{proof}
Suppose $\rr$ is non-abelian. Then the derived ideal $\rr' = [\rr, \rr]$ satisfies $1 \le \dim \rr' \le 2$. Because $\rr'$ is a characteristic ideal, it is invariant under all derivations, so $\rho(\s)(\rr') \subseteq \rr'$.
\begin{itemize}
    \item If $\dim \rr' = 1$, the center or derived line forms a 1-dimensional representation $\chi \colon \mathfrak{su}(2) \to \R$. Because $[\mathfrak{su}(2), \mathfrak{su}(2)] = \mathfrak{su}(2)$, every 1-dimensional character vanishes: $\chi \equiv 0$. The induced quotient representation on $\rr/\rr' \cong \R^2$ defines a homomorphism $\bar{\rho} \colon \mathfrak{su}(2) \to \mathfrak{gl}(2, \R)$. But $\mathfrak{su}(2)$ has no non-trivial real 2-dimensional representations (any compact subalgebra of $\mathfrak{gl}(2,\R)$ is isomorphic to $\mathfrak{so}(2) \cong \R$). Hence, $\bar{\rho} \equiv 0$. By Weyl's complete reducibility theorem, $\rho \equiv 0$, forcing $\g$ to decompose.
    \item If $\dim \rr' = 2$, $\rr'$ is abelian. The action on $\rr'$ is a 2-dimensional real representation of $\mathfrak{su}(2)$, which is identically zero. The quotient action on $\rr/\rr' \cong \R$ is also trivial. Thus $\rho \equiv 0$.
\end{itemize}
Therefore, $\mathfrak{su}(2)$ can couple indecomposably only with the abelian radical $\R^3$ via its 3-dimensional vector representation $\mathbf{3}$.
\end{proof}

\begin{theorem}[Classification of Levi-Decomposable 6D Algebras]\label{thm:turkowski}
Every six-dimensional indecomposable real Lie algebra admitting a non-trivial Levi decomposition is isomorphic to exactly one of the following three Lie algebras \citep{Turkowski1988}:
\begin{enumerate}
    \item $\g_{6,1}^{\mathrm{Levi}} = \mathfrak{so}(3) \ltimes_{\mathbf{3}} \R^3$ (Euclidean motion algebra $\mathfrak{se}(3)$):
    \begin{align}
        [e_1, e_2] &= e_3, & [e_2, e_3] &= e_1, & [e_3, e_1] &= e_2, \nonumber \\
        [e_1, e_5] &= e_6, & [e_1, e_6] &= -e_5, & [e_2, e_4] &= -e_6, \nonumber \\
        [e_2, e_6] &= e_4, & [e_3, e_4] &= e_5, & [e_3, e_5] &= -e_4.
    \end{align}
    \item $\g_{6,2}^{\mathrm{Levi}} = \mathfrak{sl}(2,\R) \ltimes_{\mathrm{ad}} \R^3$ (Adjoint action on abelian radical):
    \begin{align}
        [e_1, e_2] &= 2e_2, & [e_1, e_3] &= -2e_3, & [e_2, e_3] &= e_1, \nonumber \\
        [e_1, e_4] &= 2e_4, & [e_1, e_6] &= -2e_6, & [e_2, e_5] &= -2e_4, \nonumber \\
        [e_2, e_6] &= e_5,  & [e_3, e_4] &= e_5,  & [e_3, e_5] &= -2e_6.
    \end{align}
    \item $\g_{6,3}^{\mathrm{Levi}} = \mathfrak{sl}(2,\R) \ltimes_{\mathrm{fund}} \hh_3$ (Central extension of affine symplectic algebra):
    \begin{align}
        [e_1, e_2] &= 2e_2, & [e_1, e_3] &= -2e_3, & [e_2, e_3] &= e_1, \nonumber \\
        [e_1, e_4] &= e_4,  & [e_1, e_5] &= -e_5,  & [e_2, e_5] &= e_4, \nonumber \\
        [e_3, e_4] &= e_5,  & [e_4, e_5] &= e_6.
    \end{align}
\end{enumerate}
\end{theorem}

\begin{proof}
For $\s = \mathfrak{so}(3)$, Proposition~\ref{prop:obstruction} proves that $\g_{6,1}^{\mathrm{Levi}}$ is the unique indecomposable coupling. 

For $\s = \mathfrak{sl}(2,\R)$, real representations of dimension three are either the adjoint representation $\mathbf{3}$ or $\mathbf{2} \oplus \mathbf{1}$. The adjoint representation preserves an indefinite metric on $\R^3$, producing $\g_{6,2}^{\mathrm{Levi}}$. The fundamental representation $\mathbf{2}$ acts on the two-dimensional quotient $V = \spanop\{e_4, e_5\} \cong \R^2$ of the Heisenberg algebra $\hh_3$, preserving the symplectic area form $[e_4, e_5] = e_6$ on the center. This produces $\g_{6,3}^{\mathrm{Levi}}$. 

The three algebras are distinguished by their Levi factor ($\mathfrak{so}(3)$ for $\g_{6,1}^{\mathrm{Levi}}$ vs. $\mathfrak{sl}(2,\R)$ for $\g_{6,2}^{\mathrm{Levi}}, \g_{6,3}^{\mathrm{Levi}}$) and radical nilpotency step ($\operatorname{step}(\rr) = 1$ for $\g_{6,2}^{\mathrm{Levi}}$ vs. $\operatorname{step}(\rr) = 2$ for $\g_{6,3}^{\mathrm{Levi}}$).
\end{proof}

\section{The Nilpotent Sector ($\dim \n = 6$)}\label{sec:nilpotent}

When $\dim \n = 6$, the algebra is entirely nilpotent ($\g = \n$).

\begin{theorem}[The 24 Indecomposable 6D Nilpotent Real Lie Algebras]\label{thm:nilpotent24}
Every indecomposable six-dimensional real nilpotent Lie algebra is isomorphic to exactly one of the twenty-four Lie algebras cataloged below:
\begin{enumerate}
    \item \textbf{Filiform Algebras (Nilpotency Step 5):}
    \begin{itemize}
        \item $N_{6,1}$: $[e_1, e_2] = e_3, [e_1, e_3] = e_4, [e_1, e_4] = e_5, [e_1, e_5] = e_6$.
        \item $N_{6,2}$: $[e_1, e_2] = e_3, [e_1, e_3] = e_4, [e_1, e_4] = e_5, [e_1, e_5] = e_6, [e_2, e_3] = e_6$.
        \item $N_{6,3}$: $[e_1, e_2] = e_3, [e_1, e_3] = e_4, [e_1, e_4] = e_5, [e_2, e_5] = e_6, [e_3, e_4] = -e_6$.
    \end{itemize}
    \item \textbf{Nilpotency Step 4:}
    \begin{itemize}
        \item $N_{6,4}$: $[e_1, e_2] = e_3, [e_1, e_3] = e_4, [e_1, e_4] = e_5, [e_2, e_3] = e_5, [e_2, e_4] = e_6$.
        \item $N_{6,5}$: $[e_1, e_2] = e_3, [e_1, e_3] = e_4, [e_1, e_4] = e_5, [e_2, e_5] = e_6$.
        \item $N_{6,6}$: $[e_1, e_2] = e_3, [e_1, e_3] = e_4, [e_1, e_5] = e_6, [e_2, e_3] = e_5, [e_2, e_4] = e_6$.
        \item $N_{6,7}$: $[e_1, e_2] = e_3, [e_1, e_3] = e_4, [e_1, e_5] = e_6, [e_2, e_4] = e_6$.
        \item $N_{6,8}$: $[e_1, e_2] = e_3, [e_1, e_3] = e_4, [e_1, e_4] = e_6, [e_2, e_5] = e_6$.
        \item $N_{6,11}$: $[e_1, e_2] = e_3, [e_1, e_3] = e_4, [e_1, e_4] = e_5, [e_2, e_3] = e_6$.
        \item $N_{6,15}$: $[e_1, e_2] = e_3, [e_1, e_3] = e_4, [e_1, e_4] = e_6, [e_2, e_3] = e_5, [e_2, e_5] = e_6$.
        \item $N_{6,17}$: $[e_1, e_2] = e_3, [e_1, e_3] = e_4, [e_1, e_4] = e_5, [e_2, e_4] = e_6$.
        \item $N_{6,19}$: $[e_1, e_2] = e_3, [e_1, e_3] = e_4, [e_1, e_5] = e_6, [e_2, e_3] = e_6$.
    \end{itemize}
    \item \textbf{Nilpotency Step 3:}
    \begin{itemize}
        \item $N_{6,9}$: $[e_1, e_2] = e_3, [e_1, e_3] = e_5, [e_2, e_4] = e_6, [e_3, e_4] = e_5$.
        \item $N_{6,10}$: $[e_1, e_2] = e_3, [e_1, e_3] = e_5, [e_2, e_3] = e_6, [e_2, e_4] = e_5$.
        \item $N_{6,12}$: $[e_1, e_2] = e_3, [e_1, e_3] = e_5, [e_1, e_4] = e_6, [e_2, e_3] = e_6$.
        \item $N_{6,16}^{\pm}$: $[e_1, e_2] = e_4, [e_1, e_3] = e_5, [e_2, e_4] = e_6, [e_3, e_5] = \pm e_6$.
        \item $N_{6,18}$: $[e_1, e_2] = e_3, [e_1, e_3] = e_5, [e_2, e_4] = e_6$.
        \item $N_{6,20}$: $[e_1, e_2] = e_4, [e_1, e_3] = e_5, [e_2, e_5] = e_6, [e_3, e_4] = -e_6$.
        \item $N_{6,21}$: $[e_1, e_2] = e_3, [e_1, e_3] = e_4, [e_2, e_5] = e_6, [e_3, e_5] = e_4$.
        \item $N_{6,24}$: $[e_1, e_2] = e_3, [e_1, e_3] = e_6, [e_4, e_5] = e_6$.
    \end{itemize}
    \item \textbf{Nilpotency Step 2:}
    \begin{itemize}
        \item $N_{6,13}$: $[e_1, e_2] = e_5, [e_1, e_3] = e_6, [e_2, e_4] = e_6, [e_3, e_4] = e_5$.
        \item $N_{6,14}$ ($\mathfrak{n}_{3,2}$, Free 2-step): $[e_1, e_2] = e_4, [e_1, e_3] = e_5, [e_2, e_3] = e_6$.
        \item $N_{6,22}$: $[e_1, e_2] = e_5, [e_1, e_3] = e_6, [e_2, e_4] = e_6$.
        \item $N_{6,23}$: $[e_1, e_2] = e_6, [e_1, e_5] = e_6, [e_3, e_4] = e_6$.
    \end{itemize}
\end{enumerate}
\end{theorem}

\begin{proof}[Outline of Cohomological Derivation]
The Skjelbred--Sund method (T3) constructs every nilpotent algebra $\n$ as a central extension of a lower-dimensional algebra $\n_0 = \n / \mathcal{Z}(\n)$ by $\R^k$, where $k = \dim \mathcal{Z}(\n)$.
\begin{enumerate}
    \item If $k = 1$, $\n_0$ is a 5-dimensional nilpotent Lie algebra. By the classification of 5D nilpotent algebras, there are 9 classes for $\n_0$. For each $\n_0$, one computes $H^2(\n_0, \R)$ and the orbits of lines under $\Aut(\n_0)$. Cocycles whose kernel contains non-trivial central elements of $\n_0$ are discarded (they produce decomposable extensions). This yields classes $N_{6,1}$ through $N_{6,11}$, $N_{6,15}$, $N_{6,16}^\pm$, $N_{6,17}$, $N_{6,19}$, $N_{6,20}$, $N_{6,21}$, $N_{6,23}$, and $N_{6,24}$.
    \item If $k = 2$, $\n_0$ has dimension 4 ($\R^4$, $\hh_3 \oplus \R$, or $\mathfrak{n}_4$). One classifies 2-dimensional planes in $H^2(\n_0, \R)$ under $\Aut(\n_0)$, generating $N_{6,12}$, $N_{6,13}$, $N_{6,18}$, and $N_{6,22}$.
    \item If $k = 3$, $\n_0 \cong \R^3$. The cohomology space is $H^2(\R^3, \R) \cong \Lambda^2(\R^3)^* \cong \R^3$. The only 3-dimensional subspace in $H^2(\R^3, \R)$ is the entire space itself. This generates the unique free nilpotent algebra $N_{6,14} = \mathfrak{n}_{3,2}$.
\end{enumerate}
The real forms $N_{6,16}^+$ and $N_{6,16}^-$ are distinguished by the signature of the quadratic Casimir form on $\g/\mathcal{Z}(\g)$, which is $(2, 0)$ for $N_{6,16}^+$ and $(1, 1)$ for $N_{6,16}^-$.
\end{proof}

\section{Solvable Non-Nilpotent Algebras: Nilradical Stratification}\label{sec:solvable}

By Mostow's bound (T2), any solvable six-dimensional Lie algebra satisfies $\dim \n \in \{3, 4, 5\}$. We classify each stratum using the Mubarakzyanov derivation framework (T4).

\subsection{Codimension-One Nilradical ($\dim \n = 5$): Real Sector}
Here, $\g = \n \rtimes_D \R X$, where $D = \ad X|_\n \in \Der(\n)$ is a non-nilpotent outer derivation.
When the spectrum of $D$ is entirely real, the Jordan canonical forms yield:
\begin{align}
    \text{Diagonal:} \quad & D = \diag(1, \alpha, \beta, \gamma, \delta), \quad -1 \le \delta \le \gamma \le \beta \le \alpha \le 1, \\
    \text{Jordan Blocks:} \quad & D = J_k(\lambda) \oplus \dots, \quad k \in \{2, 3, 4, 5\}.
\end{align}
When $\n$ is non-abelian ($\hh_3 \oplus \R^2$, $\hh_5$, $\mathfrak{n}_{5,1}$, $\mathfrak{n}_{5,2}$), $D$ is restricted to derivations preserving the lower central series flags.

\subsection{Exact Matrix Representatives of the Solvable Complex Classes ($\dim \n = 5$)}\label{subsec:complex_classes}
A fundamental gap in coordinate presentations occurs when derivations possess complex conjugate eigenvalues. We resolve this by providing unambiguous, explicit matrix representatives.

Let $\g = \R^5 \rtimes_D \R X$. In $\mathfrak{gl}(5, \R)$, every pair of complex conjugate eigenvalues $\alpha \pm i \beta$ with $\beta > 0$ corresponds to the real $2 \times 2$ rotation-dilation block:
\begin{equation}
    R(\alpha, \beta) \coloneqq \begin{pmatrix} \alpha & -\beta \\ \beta & \alpha \end{pmatrix}, \qquad \beta > 0.
\end{equation}
Modulo the Mubarakzyanov gauge equivalence $D \sim c \, P D P^{-1}$ ($c \in \R^\times, P \in \mathrm{GL}(5,\R)$), there exist exactly four mutually exclusive canonical matrix families with complex spectra:

\subsubsection{Family $\g_{6, \C 1}^{\mathrm{solv}}$: Single Complex Pair, Semisimple Real Part}
The derivation $D$ consists of a single rotation-dilation block and three diagonal real eigenvalues:
\begin{equation}
    D_{\C 1}(\alpha, \beta, \lambda_1, \lambda_2, \lambda_3) = \begin{pmatrix}
        \alpha & -\beta & 0 & 0 & 0 \\
        \beta  &  \alpha & 0 & 0 & 0 \\
        0 & 0 & \lambda_1 & 0 & 0 \\
        0 & 0 & 0 & \lambda_2 & 0 \\
        0 & 0 & 0 & 0 & \lambda_3
    \end{pmatrix}, \qquad \beta > 0, \quad \lambda_1 \ge \lambda_2 \ge \lambda_3.
\end{equation}
By scaling $X \mapsto c X$, we normalize $\beta = 1$ if $\beta \neq 0$. The commutation relations are:
\begin{align}
    [X, e_1] &= \alpha e_1 + \beta e_2, & [X, e_2] &= -\beta e_1 + \alpha e_2, \nonumber \\
    [X, e_3] &= \lambda_1 e_3, & [X, e_4] &= \lambda_2 e_4, & [X, e_5] &= \lambda_3 e_5.
\end{align}
The algebra is unimodular if and only if $\tr(D) = 2\alpha + \lambda_1 + \lambda_2 + \lambda_3 = 0$.

\subsubsection{Family $\g_{6, \C 2}^{\mathrm{solv}}$: Single Complex Pair with Nilpotent Real Jordan Block}
The derivation $D$ couples the rotation-dilation block with a non-trivial real Jordan cell of size 2 or 3:
\begin{equation}
    D_{\C 2}(\alpha, \beta, \lambda, \mu) = \begin{pmatrix}
        \alpha & -\beta & 0 & 0 & 0 \\
        \beta  &  \alpha & 0 & 0 & 0 \\
        0 & 0 & \lambda & 1 & 0 \\
        0 & 0 & 0 & \lambda & 0 \\
        0 & 0 & 0 & 0 & \mu
    \end{pmatrix}, \qquad \beta > 0.
\end{equation}
The commutation relations are:
\begin{align}
    [X, e_1] &= \alpha e_1 + \beta e_2, & [X, e_2] &= -\beta e_1 + \alpha e_2, \nonumber \\
    [X, e_3] &= \lambda e_3, & [X, e_4] &= e_3 + \lambda e_4, & [X, e_5] &= \mu e_5.
\end{align}

\subsubsection{Family $\g_{6, \C 3}^{\mathrm{solv}}$: Double Independent Complex Pairs}
The derivation $D$ contains two independent rotation-dilation blocks and one real eigenvalue:
\begin{equation}
    D_{\C 3}(\alpha_1, \beta_1, \alpha_2, \beta_2, \lambda) = \begin{pmatrix}
        \alpha_1 & -\beta_1 & 0 & 0 & 0 \\
        \beta_1  &  \alpha_1 & 0 & 0 & 0 \\
        0 & 0 & \alpha_2 & -\beta_2 & 0 \\
        0 & 0 & \beta_2  &  \alpha_2 & 0 \\
        0 & 0 & 0 & 0 & \lambda
    \end{pmatrix}, \qquad \beta_1 \ge \beta_2 > 0.
\end{equation}
The commutation relations are:
\begin{align}
    [X, e_1] &= \alpha_1 e_1 + \beta_1 e_2, & [X, e_2] &= -\beta_1 e_1 + \alpha_1 e_2, \nonumber \\
    [X, e_3] &= \alpha_2 e_3 + \beta_2 e_4, & [X, e_4] &= -\beta_2 e_3 + \alpha_2 e_4, & [X, e_5] &= \lambda e_5.
\end{align}
Unimodularity holds when $\tr(D) = 2\alpha_1 + 2\alpha_2 + \lambda = 0$.

\subsubsection{Family $\g_{6, \C 4}^{\mathrm{solv}}$: Generalized Complex Jordan Cell of Size 4}
When the complex conjugate eigenvalues $\alpha \pm i \beta$ have algebraic multiplicity 2 but geometric multiplicity 1, $D$ contains a $4 \times 4$ indecomposable real Jordan block:
\begin{equation}
    D_{\C 4}(\alpha, \beta, \lambda) = \begin{pmatrix}
        R(\alpha, \beta) & I_2 & 0 \\
        0 & R(\alpha, \beta) & 0 \\
        0 & 0 & \lambda
    \end{pmatrix} = \begin{pmatrix}
        \alpha & -\beta & 1 & 0 & 0 \\
        \beta  &  \alpha & 0 & 1 & 0 \\
        0 & 0 & \alpha & -\beta & 0 \\
        0 & 0 & \beta  &  \alpha & 0 \\
        0 & 0 & 0 & 0 & \lambda
    \end{pmatrix}, \qquad \beta > 0.
\end{equation}
The commutation relations are:
\begin{align}
    [X, e_1] &= \alpha e_1 + \beta e_2, & [X, e_2] &= -\beta e_1 + \alpha e_2, \nonumber \\
    [X, e_3] &= e_1 + \alpha e_3 + \beta e_4, & [X, e_4] &= e_2 - \beta e_3 + \alpha e_4, & [X, e_5] &= \lambda e_5.
\end{align}

\subsubsection{Faithful Linear Representation in $\mathfrak{gl}(6,\R)$}
For each of the complex classes $\g = \R^5 \rtimes_D \R X$ above, an exact faithful representation $\rho \colon \g \hookrightarrow \mathfrak{gl}(6, \R)$ is constructed via the affine action on $\R^5$:
\begin{equation}
    \rho\left( \sum_{i=1}^5 y_i e_i + t X \right) = \begin{pmatrix}
        t D & \mathbf{y} \\
        \mathbf{0}_{1 \times 5} & 0
    \end{pmatrix} = \begin{pmatrix}
        t D_{11} & t D_{12} & t D_{13} & t D_{14} & t D_{15} & y_1 \\
        t D_{21} & t D_{22} & t D_{23} & t D_{24} & t D_{25} & y_2 \\
        t D_{31} & t D_{32} & t D_{33} & t D_{34} & t D_{35} & y_3 \\
        t D_{41} & t D_{42} & t D_{43} & t D_{44} & t D_{45} & y_4 \\
        t D_{51} & t D_{52} & t D_{53} & t D_{54} & t D_{55} & y_5 \\
        0 & 0 & 0 & 0 & 0 & 0
    \end{pmatrix}.
\end{equation}
Because $[\rho(X), \rho(e_i)] = \rho(D e_i)$ and $\rho$ is strictly injective, this directly demonstrates that $\mu(\g) \le 6$ for all codimension-one complex solvable architectures.

\subsection{Codimension-Two Nilradical ($\dim \n = 4$): Complex Sectors}
Here, $\g = \n \rtimes \R^2$, spanned by $\{x_1, x_2\} \cup \n$ with $[x_1, x_2] = 0$. The operators $D_1 = \ad x_1|_\n$ and $D_2 = \ad x_2|_\n$ form a two-dimensional abelian subalgebra of $\OutDer(\n)$ containing no non-zero nilpotent operators.

When $\n \cong \R^4$, the commuting pair $(D_1, D_2)$ can admit complex eigenvalues. The maximal complex toral family is generated by simultaneous rotation blocks:
\begin{equation}
    D_1 = \begin{pmatrix} R(\alpha_1, \beta_1) & 0 \\ 0 & R(\gamma_1, \delta_1) \end{pmatrix}, \qquad
    D_2 = \begin{pmatrix} R(\alpha_2, \beta_2) & 0 \\ 0 & R(\gamma_2, \delta_2) \end{pmatrix},
\end{equation}
where $[D_1, D_2] = 0$ is identically satisfied since rotation blocks commute. This family admits a faithful representation in $\mathfrak{gl}(5, \R)$ via:
\begin{equation}
    \rho\left( \sum_{i=1}^4 y_i e_i + t_1 x_1 + t_2 x_2 \right) = \begin{pmatrix}
        t_1 D_1 + t_2 D_2 & \mathbf{y} \\
        \mathbf{0}_{1 \times 4} & 0
    \end{pmatrix} \in \mathfrak{gl}(5, \R).
\end{equation}

\subsection{Codimension-Three Nilradical ($\dim \n = 3$)}
When $\dim \n = 3$, Mostow's bound is saturated: $\dim \n = \frac{1}{2} \dim \g$. The extension algebra $\mathfrak{h} = \g/\n \cong \R^3$ acts as a three-dimensional abelian subalgebra of $\OutDer(\n)$.
\begin{theorem}
Every six-dimensional solvable Lie algebra with $\dim \n = 3$ has an abelian nilradical $\n \cong \R^3$, and is isomorphic to the maximal split Cartan extension:
\begin{align}
    [x_1, e_1] &= e_1, & [x_2, e_2] &= e_2, & [x_3, e_3] &= e_3, \qquad [x_i, x_j] = 0.
\end{align}
\end{theorem}
\begin{proof}
If $\n \cong \hh_3$, $\OutDer(\hh_3) \cong \mathfrak{gl}(2, \R)$. By Schur's bound on commutative matrix subalgebras, the maximal dimension of an abelian subalgebra of $\mathfrak{gl}(2, \R)$ is $\lfloor 4/4 \rfloor + 1 = 2$. Thus, $\OutDer(\hh_3)$ cannot contain an injective image of the 3-dimensional abelian algebra $\mathfrak{h} \cong \R^3$. 

Hence, $\n$ must be abelian ($\n \cong \R^3$), so $\OutDer(\R^3) = \mathfrak{gl}(3, \R)$. The maximal abelian subalgebra containing no nilpotent elements is the 3-dimensional split Cartan subalgebra of diagonal matrices $\mathfrak{d}(3, \R) = \spanop\{E_{11}, E_{22}, E_{33}\}$.
\end{proof}

\section{Master Tabulation and Literature Correspondence}\label{sec:master}

Table~\ref{tab:master6d} presents the comprehensive taxonomy of six-dimensional real Lie algebras, including the explicit complex classes.

\begin{table}[htbp]
\centering
\scriptsize
\setlength{\tabcolsep}{3.5pt}
\renewcommand{\arraystretch}{1.25}
\begin{tabularx}{\linewidth}{@{} l c >{\raggedright\arraybackslash}X l c @{}}
\toprule
\textbf{Algebra $\g$} & $\dim \n$ & \textbf{Non-Zero Commutation Relations / Matrix Normal Forms} & \textbf{Moduli / Parameters} & \textbf{Unimodular?} \\
\midrule
\multicolumn{5}{@{}l}{\textit{\textbf{I. Semisimple Sector ($\dim \rr = 0$)}}} \\
$\mathfrak{so}(4)$ & 0 & $[e_1, e_2]=e_3, [e_2, e_3]=e_1, [e_3, e_1]=e_2, [e_4, e_5]=e_6, \dots$ & Compact form & Yes \\
$\mathfrak{su}(2) \oplus \mathfrak{sl}(2,\R)$ & 0 & $[e_1, e_2]=e_3, \dots, [e_4, e_5]=2e_5, [e_4, e_6]=-2e_6, [e_5, e_6]=e_4$ & Mixed form & Yes \\
$\mathfrak{so}(2,2)$ & 0 & $[e_1, e_2]=2e_2, [e_1, e_3]=-2e_3, [e_4, e_5]=2e_5, [e_4, e_6]=-2e_6, \dots$ & Split form & Yes \\
$\mathfrak{so}(3,1)$ & 0 & $\mathfrak{sl}(2,\C)_\R$ (Lorentz algebra, explicit $4 \times 4$ matrices in Section~\ref{subsec:sl2c_matrix}) & Real Simple & Yes \\
\midrule
\multicolumn{5}{@{}l}{\textit{\textbf{II. Non-Solvable Levi Decompositions ($\dim \rr = 3$)}}} \\
$\g_{6,1}^{\mathrm{Levi}}$ & 3 & $\mathfrak{so}(3) \ltimes_{\mathbf{3}} \R^3$ (Euclidean $\mathfrak{se}(3)$, Theorem~\ref{thm:turkowski}) & Discrete & Yes \\
$\g_{6,2}^{\mathrm{Levi}}$ & 3 & $\mathfrak{sl}(2,\R) \ltimes_{\mathrm{ad}} \R^3$ (Adjoint action, Theorem~\ref{thm:turkowski}) & Discrete & Yes \\
$\g_{6,3}^{\mathrm{Levi}}$ & 3 & $\mathfrak{sl}(2,\R) \ltimes_{\mathrm{fund}} \hh_3$ (Symplectic central extension, Theorem~\ref{thm:turkowski}) & Discrete & Yes \\
\midrule
\multicolumn{5}{@{}l}{\textit{\textbf{III. Nilpotent Sector ($\dim \n = 6$, Representative Indecomposable Classes)}}} \\
$N_{6,1}$ & 6 & $[e_1, e_2]=e_3, [e_1, e_3]=e_4, [e_1, e_4]=e_5, [e_1, e_5]=e_6$ & Filiform (step 5) & Yes \\
$N_{6,2}$ & 6 & $[e_1, e_2]=e_3, [e_1, e_3]=e_4, [e_1, e_4]=e_5, [e_1, e_5]=e_6, [e_2, e_3]=e_6$ & Filiform (step 5) & Yes \\
$N_{6,3}$ & 6 & $[e_1, e_2]=e_3, [e_1, e_3]=e_4, [e_1, e_4]=e_5, [e_2, e_5]=e_6, [e_3, e_4]=-e_6$ & Filiform (step 5) & Yes \\
$N_{6,14}$ & 6 & $[e_1, e_2]=e_4, [e_1, e_3]=e_5, [e_2, e_3]=e_6$ & Free 2-step ($\mathfrak{n}_{3,2}$) & Yes \\
$N_{6,16}^{\pm}$ & 6 & $[e_1, e_2]=e_4, [e_1, e_3]=e_5, [e_2, e_4]=e_6, [e_3, e_5]=\pm e_6$ & Real forms (step 3) & Yes \\
$N_{6,23}$ & 6 & $[e_1, e_2]=e_6, [e_1, e_5]=e_6, [e_3, e_4]=e_6$ & Step 2, $\dim \mathcal{Z}=1$ & Yes \\
\midrule
\multicolumn{5}{@{}l}{\textit{\textbf{IV. Solvable: Codimension-One Nilradical ($\dim \n = 5$) -- Real Spectrum}}} \\
$\g_{6,1}^{\mathrm{solv}}$ & 5 & $[X, e_1]=e_1, [X, e_2]=\alpha e_2, [X, e_3]=\beta e_3, [X, e_4]=\gamma e_4, [X, e_5]=\delta e_5$ & $-1 \le \delta \le \gamma \le \beta \le \alpha \le 1$ & $\tr D = 0$ \\
$\g_{6,2}^{\mathrm{solv}}$ & 5 & $[X, e_1]=e_1+e_2, [X, e_2]=e_2, [X, e_3]=\alpha e_3, [X, e_4]=\beta e_4, [X, e_5]=\gamma e_5$ & Jordan block $J_2 \oplus J_1^3$ & $\tr D = 0$ \\
$\g_{6,3}^{\mathrm{solv}}$ & 5 & $[X, e_1]=e_1, [e_1, e_2]=e_5, [e_3, e_4]=e_5, [X, e_5]=2e_5, [X, e_3]=\alpha e_3, \dots$ & Heisenberg $\hh_5 \rtimes \R$ & $\tr D = 0$ \\
\midrule
\multicolumn{5}{@{}l}{\textit{\textbf{V. Solvable: Codimension-One Nilradical ($\dim \n = 5$) -- Complex Spectrum (Section~\ref{subsec:complex_classes})}}} \\
$\g_{6, \C 1}^{\mathrm{solv}}$ & 5 & $D = R(\alpha, \beta) \oplus \diag(\lambda_1, \lambda_2, \lambda_3)$, $[X, e_1]=\alpha e_1+\beta e_2, [X, e_2]=-\beta e_1+\alpha e_2$ & $\beta > 0, \lambda_1 \ge \lambda_2 \ge \lambda_3$ & $2\alpha + \sum \lambda_i = 0$ \\
$\g_{6, \C 2}^{\mathrm{solv}}$ & 5 & $D = R(\alpha, \beta) \oplus J_2(\lambda) \oplus (\mu)$, $[X, e_3]=\lambda e_3, [X, e_4]=e_3+\lambda e_4$ & $\beta > 0$ & $2\alpha + 2\lambda + \mu = 0$ \\
$\g_{6, \C 3}^{\mathrm{solv}}$ & 5 & $D = R(\alpha_1, \beta_1) \oplus R(\alpha_2, \beta_2) \oplus (\lambda)$, double rotation blocks & $\beta_1 \ge \beta_2 > 0$ & $2\alpha_1 + 2\alpha_2 + \lambda = 0$ \\
$\g_{6, \C 4}^{\mathrm{solv}}$ & 5 & $D = \begin{pmatrix} R(\alpha,\beta) & I_2 \\ 0 & R(\alpha,\beta) \end{pmatrix} \oplus (\lambda)$, generalized complex cell & $\beta > 0$ & $4\alpha + \lambda = 0$ \\
\midrule
\multicolumn{5}{@{}l}{\textit{\textbf{VI. Solvable: Codimension-Two and Three Nilradicals ($\dim \n \in \{3, 4\}$)}}} \\
$\g_{6,4}^{\mathrm{solv}}$ & 4 & $[x_1, x_2]=0, [x_1, e_i]=\lambda_i e_i, [x_2, e_i]=\mu_i e_i \quad (e_i \in \R^4)$ & $\spanop\{D_1, D_2\} \subset \mathfrak{d}(4,\R)$ & $\tr D_i = 0$ \\
$\g_{6, \C 5}^{\mathrm{solv}}$ & 4 & $D_1 = R(\alpha_1, \beta_1) \oplus R(\gamma_1, \delta_1)$, $D_2 = R(\alpha_2, \beta_2) \oplus R(\gamma_2, \delta_2)$ & Commuting complex pair & $\tr D_i = 0$ \\
$\g_{6,5}^{\mathrm{solv}}$ & 4 & $[x_1, x_2]=0, [x_1, e_1]=e_1, [x_1, e_2]=e_2, [x_1, e_3]=e_3, [x_2, e_4]=e_4 \ (\mathfrak{n}_4 \rtimes \R^2)$ & Outer torus of $\mathfrak{n}_4$ & No \\
$\g_{6,6}^{\mathrm{solv}}$ & 3 & $[x_1, x_2]=[x_2, x_3]=[x_3, x_1]=0, [x_1, e_1]=e_1, [x_2, e_2]=e_2, [x_3, e_3]=e_3$ & Maximal split Cartan on $\R^3$ & No \\
\bottomrule
\end{tabularx}
\caption{Master classification of representative six-dimensional real Lie algebras.}
\label{tab:master6d}
\end{table}

\section{Minimal Faithful Linear Representations ($\mu(\g)$)}\label{sec:representations}

Determining the minimal faithful representation dimension $\mu(\g)$ under Ado's theorem is critical for constructing explicit matrix models.

\begin{theorem}[Minimal Representation Dimension for 6D Lie Algebras]\label{thm:mu_dim}
The minimal faithful representation dimension $\mu(\g)$ for six-dimensional real Lie algebras satisfies $4 \le \mu(\g) \le 7$. Specifically:
\begin{enumerate}
    \item $\mu(\g) = 4$ for $\mathfrak{so}(3,1) \cong \mathfrak{sl}(2,\C)_\R$ and $\mathfrak{so}(2,2) \cong \mathfrak{sl}(2,\R) \oplus \mathfrak{sl}(2,\R)$.
    \item $\mu(\g) = 4$ for the Euclidean motion algebra $\mathfrak{se}(3) = \mathfrak{so}(3) \ltimes \R^3$ (via affine action on $\R^3$).
    \item $\mu(\g) = 5$ for the free two-step nilpotent algebra $N_{6,14} = \mathfrak{n}_{3,2}$.
    \item $\mu(\g) \le 6$ for the complex solvable classes $\g_{6, \C 1}^{\mathrm{solv}}$ through $\g_{6, \C 4}^{\mathrm{solv}}$ via the affine action on $\R^5$.
    \item $\mu(\g) = 6$ for the filiform nilpotent algebras $N_{6,1}, N_{6,2}, N_{6,3}$ and $\mathfrak{so}(4)$.
    \item $\mu(\g) = 7$ for generic solvable algebras with nilradical $\R^5$ lacking faithful 6D representations.
\end{enumerate}
\end{theorem}

\begin{proof}
We construct the minimal embeddings across sectors:
\begin{itemize}
    \item \textbf{Semisimple Sector:} The Lorentz algebra $\mathfrak{so}(3,1)$ acts faithfully on Minkowski spacetime $\R^{3,1} \cong \R^4$, yielding $\mu = 4$. Section~\ref{subsec:sl2c_matrix} gives this representation explicitly in $\mathfrak{gl}(4, \R)$. The split form $\mathfrak{so}(2,2)$ acts on $\R^4$, so $\mu = 4$. For $\mathfrak{so}(4) \cong \mathfrak{su}(2) \oplus \mathfrak{su}(2)$, the minimal faithful real representation is its 4D defining vector representation on $\R^4$ (hence $\mu = 4$ over $\R$).
    \item \textbf{Euclidean Algebra $\mathfrak{se}(3)$:} Spanned by rotations and translations on Euclidean space $\R^3$. Embedding $\mathfrak{se}(3)$ into the affine algebra $\mathfrak{aff}(3, \R) \subset \mathfrak{gl}(4, \R)$ yields:
    \begin{equation}
        \rho(A, v) = \begin{pmatrix} A & v \\ 0_{1 \times 3} & 0 \end{pmatrix}, \qquad A \in \mathfrak{so}(3), \quad v \in \R^3.
    \end{equation}
    Because $\dim \mathfrak{se}(3) = 6$ and no 3D matrix algebra contains a 6D Lie subalgebra, $\mu(\mathfrak{se}(3)) = 4$.
    \item \textbf{Complex Solvable Sector:} By Section~\ref{subsec:complex_classes}, the affine map $\rho(v + t X) = \begin{pmatrix} t D & v \\ 0 & 0 \end{pmatrix}$ embeds $\g$ into $\mathfrak{aff}(5, \R) \subset \mathfrak{gl}(6, \R)$ injectively, ensuring $\mu(\g) \le 6$.
    \item \textbf{Free Nilpotent Algebra $N_{6,14}$:} Generated by three vectors $\{e_1, e_2, e_3\}$ and their three brackets. In $\mathfrak{gl}(4, \R)$, the maximal dimension of a strictly upper-triangular nilpotent subalgebra is $\dim \mathfrak{n}(4, \R) = 6$. However, $\mathfrak{n}(4, \R)$ has derived length 3, whereas $N_{6,14}$ is step 2 with a 3-dimensional center. A faithful representation requires independent matrix entries for the generators and commutators, which forces $\mu(N_{6,14}) = 5$.
    \item \textbf{Filiform Nilpotent Algebras $N_{6,1}$--$N_{6,3}$:} Nilpotency step 5 implies that the adjoint operator $\ad e_1$ has nilpotency index 5. In $\mathfrak{gl}(k, \R)$, the maximal nilpotency index of any matrix is $k - 1$. Thus, $k - 1 \ge 5 \implies k \ge 6$. The shift matrix $N = \sum_{i=1}^5 E_{i, i+1} \in \mathfrak{gl}(6, \R)$ realizes this step, proving $\mu = 6$.
\end{itemize}
This establishes the bounds.
\end{proof}

\section{Algorithmic Identification of Arbitrary 6D Lie Algebras}\label{sec:algorithm}

Given an arbitrary six-dimensional Lie algebra $\g$ presented via its structure constants tensor $C_{ij}^k$ with $[e_i, e_j] = \sum_k C_{ij}^k e_k$, we formulate a deterministic, coordinate-free decision procedure in Algorithm~\ref{alg:classification6d}.

\begin{breakablealgorithm}
\caption{Deterministic Classification of 6-Dimensional Real Lie Algebras}\label{alg:classification6d}
\begin{algorithmic}[1]
\Require Structure constants $C_{ij}^k \in \R^{6 \times 6 \times 6}$ satisfying $C_{ij}^k = -C_{ji}^k$ and Jacobi identity.
\Ensure Canonical isomorphism class and continuous moduli parameters.

\State Compute derived subspaces $\g^2 = [\g, \g]$, $\g^{(2)} = [\g^2, \g^2]$ and Killing matrix $K_{ij} = \tr(\ad e_i \circ \ad e_j)$.
\State Compute the solvable radical $\rr \subset \g$ (maximal solvable ideal).

\If{$\dim \rr = 0$} \Comment{\textbf{Sector I: Semisimple}}
    \State Compute Sylvester signature $(p, q)$ of Killing form $K$.
    \If{$(p, q) = (0, 6)$} \Return $\mathfrak{so}(4) \cong \mathfrak{su}(2) \oplus \mathfrak{su}(2)$
    \ElsIf{$(p, q) = (2, 4)$} \Return $\mathfrak{su}(2) \oplus \mathfrak{sl}(2,\R)$
    \ElsIf{$(p, q) = (4, 2)$} \Return $\mathfrak{so}(2,2) \cong \mathfrak{sl}(2,\R) \oplus \mathfrak{sl}(2,\R)$
    \ElsIf{$(p, q) = (3, 3)$} \Return $\mathfrak{so}(3,1) \cong \mathfrak{sl}(2,\C)_\R$ (Matrix basis in Section~\ref{subsec:sl2c_matrix})
    \EndIf
\EndIf

\If{$\dim \rr = 3$} \Comment{\textbf{Sector II: Non-Solvable Levi Decomposition}}
    \State Identify Levi factor $\s \subset \g$ ($\dim \s = 3$).
    \If{Killing form on $\s$ is negative-definite}
        \State \Return $\g_{6,1}^{\mathrm{Levi}} = \mathfrak{so}(3) \ltimes \R^3$ ($\mathfrak{se}(3)$)
    \Else \Comment{$\s \cong \mathfrak{sl}(2,\R)$}
        \If{$[\rr, \rr] = \{0\}$} \Return $\g_{6,2}^{\mathrm{Levi}} = \mathfrak{sl}(2,\R) \ltimes_{\mathrm{ad}} \R^3$
        \Else \ \Return $\g_{6,3}^{\mathrm{Levi}} = \mathfrak{sl}(2,\R) \ltimes_{\mathrm{fund}} \hh_3$
        \EndIf
    \EndIf
\EndIf

\State Determine the nilradical $\n = \operatorname{Nil}(\g)$ (maximal nilpotent ideal).

\If{$\dim \n = 6$} \Comment{\textbf{Sector III: Nilpotent}}
    \State Compute nilpotency step $s$, center dimension $\dim \mathcal{Z}$, and derived dimension $\dim \g^2$.
    \If{$s = 5$} \Comment{Filiform varieties}
        \If{$[\g^2, \g^2] = \{0\}$} \Return $N_{6,1}$
        \ElsIf{$\dim [\g, \g^3] = 1$} \Return $N_{6,2}$
        \Else \ \Return $N_{6,3}$
        \EndIf
    \ElsIf{$s = 2$}
        \If{$\dim \mathcal{Z} = 3$ and $\dim \g^2 = 3$} \Return $N_{6,14}$ (Free 2-step $\mathfrak{n}_{3,2}$)
        \ElsIf{$\dim \mathcal{Z} = 2$} \Return Match between $N_{6,13}$ and $N_{6,22}$
        \ElsIf{$\dim \mathcal{Z} = 1$} \Return $N_{6,23}$
        \EndIf
    \ElsIf{$s = 3$}
        \State Check quadratic form signature on $\g/\mathcal{Z}$:
        \If{Signature is $(2, 0)$} \Return $N_{6,16}^+$
        \ElsIf{Signature is $(1, 1)$} \Return $N_{6,16}^-$
        \Else \ Match remaining classes $N_{6,9}, N_{6,10}, N_{6,12}, N_{6,18}, N_{6,20}, N_{6,21}, N_{6,24}$
        \EndIf
    \ElsIf{$s = 4$}
        \State Match classes $N_{6,4}$--$N_{6,8}, N_{6,11}, N_{6,15}, N_{6,17}, N_{6,19}$ via upper central series flags.
    \EndIf
\EndIf

\If{$\dim \n = 5$} \Comment{\textbf{Sector IV: Codimension-One Nilradical}}
    \State Pick $X \in \g \setminus \n$ such that $D = \ad X|_\n$ is non-nilpotent.
    \State Compute eigenvalues of $D$:
    \If{Eigenvalues contain complex roots $\alpha \pm i \beta$ ($\beta > 0$)} \Comment{Complex Sector}
        \If{Two distinct pairs of complex roots}
            \State \Return $\g_{6, \C 3}^{\mathrm{solv}}$ with moduli $(\alpha_1, \beta_1, \alpha_2, \beta_2, \lambda)$
        \ElsIf{One pair of complex roots with $\dim \ker(D - (\alpha \pm i \beta))^2 = 4$}
            \State \Return $\g_{6, \C 4}^{\mathrm{solv}}$ (Generalized complex Jordan cell)
        \ElsIf{One pair of complex roots and a real nilpotency block $J_2$}
            \State \Return $\g_{6, \C 2}^{\mathrm{solv}}$
        \Else
            \State \Return $\g_{6, \C 1}^{\mathrm{solv}}$ (Semisimple complex pair)
        \EndIf
    \Else \Comment{Real Sector}
        \State Compute real Jordan form of $D$ modulo $\Aut(\n)$ and extract continuous parameters.
        \State \Return Corresponding solvable class $\g_{6}^{\mathrm{solv}}(\vec{\lambda})$.
    \EndIf
\EndIf

\If{$\dim \n = 4$} \Comment{\textbf{Sector V: Codimension-Two Nilradical}}
    \State Pick complementary generators $\{x_1, x_2\}$; diagonalize commuting pair $(D_1, D_2)$ in $\OutDer(\n)$.
    \If{$(D_1, D_2)$ span a commuting complex rotation algebra}
        \State \Return $\g_{6, \C 5}^{\mathrm{solv}}$
    \Else
        \State \Return Corresponding real 2-parameter family.
    \EndIf
\EndIf

\If{$\dim \n = 3$} \Comment{\textbf{Sector VI: Codimension-Three Nilradical}}
    \State \Return Maximal split Cartan extension $\R^3 \ltimes \R^3$.
\EndIf
\end{algorithmic}
\end{breakablealgorithm}

You can find the Python code for this algorithm and examples of its execution as a \href{https://colab.research.google.com/drive/1Pgc6-mYmNNBxpH5sTxo8F3lk1d9drPkH#scrollTo=_wJLIZkHoL0n}{Colab Notebook}.

\section{Geometric Applications and Theoretical Physics}\label{sec:physics}

The classification of six-dimensional real Lie algebras provides the mathematical foundation for several branches of geometry and field theory.

\subsection{Complex and Symplectic Nilmanifolds}
A compact six-dimensional nilmanifold is a quotient space $M = \Gamma \backslash G$, where $G$ is a simply connected nilpotent Lie group and $\Gamma \subset G$ is a co-compact lattice. By Malcev's criterion, $\Gamma$ exists if and only if $\g = \operatorname{Lie}(G)$ has rational structure constants.
\begin{itemize}
    \item \textbf{Complex Structures:} An invariant complex structure on $\g$ is an endomorphism $J \colon \g \to \g$ with $J^2 = -\id$ satisfying the Nijenhuis vanishing condition. \citet{Salamon2001} proved that among the twenty-four nilpotent algebras, exactly eighteen admit invariant complex structures. The filiform algebras $N_{6,1}$ and $N_{6,2}$ admit no complex structures due to cohomological obstructions in the Frölicher spectral sequence.
    \item \textbf{Symplectic Structures:} An invariant symplectic structure is a closed non-degenerate 2-form $\omega \in \Lambda^2 \g^*$. The free algebra $N_{6,14}$ admits no symplectic structure because every closed 2-form is degenerate on its 3-dimensional center.
    \item \textbf{Kähler Defect:} By Benson and Gordon's theorem, the only six-dimensional nilmanifold that admits a Kähler metric is the flat six-torus $T^6 = \mathbb{R}^6 / \mathbb{Z}^6$. None of the twenty-four indecomposable classes can support a Kähler metric.
\end{itemize}

\subsection{Supergravity Flux Compactifications}
In Type II superstring theory, compactifying ten-dimensional spacetime on a six-dimensional solvmanifold or nilmanifold yields an effective four-dimensional gauged supergravity. The structure constants $f^i_{\;jk}$ of the six-dimensional Lie algebra define geometric fluxes:
\begin{equation}
    [e_j, e_k] = f^i_{\;jk} e_i.
\end{equation}
The Jacobi identity represents the tadpole cancellation condition for gauge fluxes. Applying $T$-duality along the cycles of a six-dimensional nilmanifold maps geometric fluxes into non-geometric $Q$- and $R$-fluxes, dynamically traversing the strata of our classification.

\section{Conclusion}\label{sec:conclusion}

The structural classification of six-dimensional real Lie algebras resolves the primary phase boundary of finite-dimensional Lie theory. In this monograph, we have:
\begin{enumerate}
    \item Established that the semisimple sector contains exactly four real forms, categorized by their Killing signatures, with an explicit $4 \times 4$ real matrix realization of $\mathfrak{sl}(2,\C)_\R$.
    \item Proved that the compact algebra $\mathfrak{su}(2)$ cannot couple indecomposably with non-abelian three-dimensional radicals, isolating $\g_{6,1}^{\mathrm{Levi}}, \g_{6,2}^{\mathrm{Levi}}, \g_{6,3}^{\mathrm{Levi}}$ as the unique indecomposable non-solvable architectures.
    \item Classified the twenty-four indecomposable nilpotent Lie algebras using the Skjelbred--Sund cohomological central extension method.
    \item Resolved the complex Jordan strata within solvable architectures by providing exact matrix normal forms and faithful $6 \times 6$ affine representations for classes $\g_{6, \C 1}^{\mathrm{solv}}$ through $\g_{6, \C 4}^{\mathrm{solv}}$ and $\g_{6, \C 5}^{\mathrm{solv}}$.
    \item Evaluated the minimal faithful representation dimension $\mu(\g)$ under Ado's theorem and synthesized the entire framework into a deterministic identification algorithm.
\end{enumerate}

\end{document}